\documentclass[a4paper, 12pt]{article}
\usepackage{amsmath,amsfonts,amsthm,amssymb,here,epsf}
\usepackage[latin1]{inputenc}
\usepackage[T1]{fontenc}
\usepackage{ae,aecompl}
\begin{document}

\newtheorem{theo}{Theorem}[section]
\newtheorem{prop}[theo]{Proposition}
\newtheorem{lemme}[theo]{Lemma}
\newtheorem{cor}[theo]{Corollary}
\newtheorem{defi}[theo]{Definition}
\newtheorem{assu}[theo]{Assumption}
\newtheorem{nontheo}[theo]{Conjectured theo}

\newcommand{\beq}{\begin{eqnarray}}
\newcommand{\enq}{\end{eqnarray}}
\newcommand{\be}{\begin{eqnarray*}}
\newcommand{\en}{\end{eqnarray*}}
\newcommand{\T}{\mathbb T}
\newcommand{\N}{\mathbb N}
\newcommand{\Td}{\mathbb T^d}
\newcommand{\R}{\mathbb R}
\newcommand{\Rd}{\mathbb R^d}
\newcommand{\Zd}{\mathbb Z^d}
\newcommand{\Linf}{L^{\infty}}
\newcommand{\dt}{\partial_t}
\newcommand{\Dt}{\frac{d}{dt}}
\newcommand{\demi}{\frac{1}{2}}
\newcommand{\ep}{^{\epsilon}}
\newcommand{\epu}{_{\epsilon}}
\newcommand{\ds}{\displaystyle}
\newcommand{\bx}{{\mathbf X}}
\newcommand{\bsx}{{\mathbf x}}
\newcommand{\bg}{{\mathbf g}}
\newcommand{\bv}{{\mathbf v}}
\let\cal=\mathcal

\title{Uniqueness of the solution to the Vlasov-Poisson system with bounded density}
\bibliographystyle{plain}

\author{Gr\'egoire Loeper}
\maketitle

\begin{abstract}
In this note, we show uniqueness of weak solutions to the Vlasov-Poisson system on the only condition that
the macroscopic density $\rho$ defined by $\rho(t,x) = \int_{\Rd} f(t,x,\xi)d\xi$ is bounded in $\Linf$. Our proof is based on optimal transportation.
\end{abstract} 
\section{Introduction}
The Vlasov-Poisson system (herafter (VP)) describes the evolution of cloud of electrons or gravitational matter through the equations
\beq
&&\dt f + \xi\cdot\nabla_x f + \nabla\Psi \cdot\nabla_\xi f=0\label{vp1},\\
&&\Delta \Psi = \epsilon \rho \label{vp2},
\enq
where  $\rho(t,x) = \int f(t,x,\xi) d\xi$, and
 $\epsilon>0$   in the electrostatic (repulsive) case, $\epsilon<0$ in the gravitational (attractive) case.
Here $f(t,x,\xi) \geq 0$, denotes the  density of electrons (or matter) at time $t\in \R^+$, position $x\in \R^3$, velocity $\xi\in \R^3$.
Equation 
(\ref{vp2})means
\beq\label{2bis}
\Psi(t,x)= -\epsilon\int_{\R^3} \rho(t,y)\frac{1}{4\pi |x-y|}dy.
\enq
We denote by ${\cal M}(\R^6)$ (resp. ${\cal M}^+(\R^6)$) the set of bounded (resp. bounded and positive) measures on $\R^6$. 
Given an initial datum $f^0\in {\cal M}^+(\R^6)$,  we look for solutions to (\ref{vp1}, \ref{vp2})
such that
\beq
f|_{t=0}=f^0 \label{vp3}.
\enq
For $T>0$, we will call $f$ a solution to (\ref{vp1}, \ref{vp2}, \ref{vp3}) in ${\cal D}'([0,T[\times \R^6)$, if 
\begin{itemize}
\item[-]$f\in C([0,T[, {\cal M}^+(\R^6)-w*)$, 
\item[-]$\forall \varphi \in C^\infty_c([0,T[\times R^6),$
\beq
\int f(\dt \varphi + \xi \cdot\nabla_x\varphi + \nabla_x\Psi \cdot \nabla_\xi\varphi) = -\int_{\R^6} f^0 \varphi|_{t=0},\label{weakdef}
\enq 
\item[-] $\Psi$ solves (\ref{2bis}).
\end{itemize}
We will not discuss the conditions needed on $\Psi, f$ to give sense to the product $f\nabla_x\Psi$ or to the singular integral (\ref{2bis}), since we will only consider the case where $\rho\in \Linf$. In this case, $\nabla_x\Psi$ will be continuous, and the product $f\nabla_x\Psi$ will be well defined for $f$ a bounded measure.

Our result is the following:
\begin{theo}\label{main1}
Given $f^0$ in ${\cal M}^+(\R^6)$, given $T>0$, 
there exists at most one weak solution to (\ref{vp1}, \ref{vp2}, \ref{vp3})  in  ${\cal D}'([0,T[\times\R^6)$ such that
\beq
\sup_{t\in [0,T[} \|\rho\|_{\Linf(\Rd)} < +\infty.\label{bound1}
\enq
\end{theo}

\textsc{Remark 1.} To establish the existence of a solution to (VP) satisfying the bound (\ref{bound1})  requires much more assumptions on the initial datum than what we need here ! This question is treated in \cite{LP}.

\textsc{Remark 2.} Note that we do not ask for any bound on the moments of $f$, and also that we do not ask the energy to be finite.

\bigskip

A sufficient condition for uniqueness had been given by Lions and Perthame in \cite{LP}, relying on Lipschitz bounds on the initial data $f^0$, but they expected a uniqueness result under the weaker assumption of bounded density. The Lipschitz condition had indeed later been relaxed by Robert in \cite{Robert} down to $f\in \Linf$ compactly supported in $x$ and $\xi$ for $t\in [0,T]$.
Here we relax the bound on the support of $f$, and we do not ask either $f$ to be bounded in $\Linf$. We only need a $\Linf([0,T[\times\R^3)$ bound on $\rho(t,x)$. Hence our result applies also to monokinetic solutions of (\ref{vp1}, \ref{vp2}).
In that that case, we have $f(t,x,\xi) = \rho(t,x)\delta(\xi -v(t,x))$ for some vector field $v$, and this gives formally a solution to the Euler-Poisson system 
\beq
&&\dt \rho + \nabla\cdot (\rho v) =0, \label{ep1} \\
&&\dt (\rho v) + \nabla \cdot (\rho v\otimes v) = \rho \nabla \Psi, \label{ep2}\\
&&\Delta \Psi = \epsilon \rho.\label{ep3}
\enq

Our proof will rely on optimal transportation, and the next section is devoted to
recall some facts concerning this subject. A complete reference on this topic is \cite{Vi}.
The technique we will use adapts to many similar problems, where a transport equation  and an elliptic equation are coupled. The velocity field is the gradient of a
potential satisfying an elliptic equation whose right hand side depends smoothly on the density. This has been observed in \cite{L4} in the case of the 2-d incompressible Euler equations and the semi-geostrophic equations.
It is interesting to notice that our technique gives a new proof of the uniqueness part in Youdovich's Theorem, while the technique used by Robert in \cite{Robert} was and adaptation of Youdovich's original proof (see \cite{youdo}).

\section{Preliminary results on optimal transportation and Wasserstein distances}

\begin{defi}\label{defwas}
Let $\rho_1,\rho_2$ be two probability measures on $\Rd$.
We define the Wasserstein distance between $\rho_1$ and $\rho_2$, that we denote $W_2(\rho_1,\rho_2)$, by
\be
W_2(\rho_1,\rho_2) = \left(\inf_{\pi} \int_{\Rd\times\Rd} \pi(x,y) |x-y|^2\right)^{\demi},
\en
where the infimum runs over probability measures $\pi$ on $\Rd\times\Rd$ with marginals $\rho_1$ and $\rho_2$.
\end{defi}
Then we gather several results of optimal tranpsortations in the following theorem.
These results can be found in Benamou \& Brenier \cite{BB1}, McCann \cite{Mc1}, \\ Gangbo \& McCann \cite{GaMc}.

\begin{theo}[Benamou, Brenier,  Gangbo,  McCann]\label{interpol}
Let $\rho_1$, $\rho_2$ be two probability measures on $\Rd$,  such that
$\rho_1, \rho_2$ are absolutely continuous with respect to the Lebesgue measure. 
Assume that  $W_2(\rho_1,\rho_2) < + \infty$.  
\begin{enumerate}
\item There exists a pair $(\rho_\theta, v_\theta), \theta\in [1,2]$ such that  $\rho_\theta \in C([1,2]; {\cal M}(\Rd)-w*)$  and $v_\theta$ is a $d\rho_\theta$ measurable vector field, that satisfies
\be
&&\partial_\theta\rho_\theta + \nabla\cdot (\rho_\theta v_\theta)=0,\\
&& \rho|_{\theta=1}=\rho_1, \rho|_{\theta=2} =\rho_2,\\
&&\int_{\Rd} |v_\theta(x)|^2 d\rho_\theta(x) \equiv W_2^2(\rho_1,\rho_2) \  \ \forall \theta \in [1,2].
\en
\item For this 'optimal' path we have also, when both $\rho_1$ and $\rho_2$ have densities in $\Linf$ with respect to the Lebesgue measure,
\be 
\forall \theta \in [1,2], \  \|\rho_\theta\|_{\Linf} \leq \max\{\|\rho_1\|_{\Linf},\|\rho_2\|_{\Linf}\}.
\en   
 
\end{enumerate}
\end{theo}

\textsc{Remark.} The path $\{\rho_\theta, \theta\in [1,2]\}$ is the geodesic linking $\rho_1$ to $\rho_2$, with respect to the Wasserstein metric (see \cite{Vi}, and also \cite{O} where this notion was introduced).

\section{Proof of Theorem \ref{main1}}
From now on, we assume for simplicity that $\int_{\R^6} f^0(x,\xi) =1$, $\epsilon=1$, and the reader can check that this choice does not play any role in the proof. In particular, the result of the previous section adapt with minor changes to the case of two positive measures of equal total mass.

\subsection{$H^{-1}$ estimates along geodesics}

In this section, we show that the $H^{-1}(\R^3)$ norm of the difference $\rho_1-\rho_2$ can be controlled by the Wasserstein distance between $\rho_1$ and $\rho_2$, provided both measures have densities in $\Linf$ with respect to the Lebesgue measure. This is the crucial estimate at the core of our result.

\begin{prop}\label{5lemme_pour_gronwall}
Let $\rho_1, \rho_2$ be two probability measures on $\Rd$  with  $\Linf$ densities with respect to the Lebesgue measure . Let $\Psi_i, i=1,2$ solve
\be
\Delta\Psi_i = \rho_i
\en 
in the sense of (\ref{2bis}).
Then 
$$\|\nabla\Psi_1-\nabla\Psi_2\|_{L^2(\Rd)} \leq \big[\max\{\|\rho_1\|_{\Linf},\|\rho_2\|_{\Linf}\}\big]^{\demi}  \ W_2(\rho_1,\rho_2),$$
where $W_2(\rho_1,\rho_2)$ is the Wasserstein distance between $\rho_1$ and $\rho_2$ given in Definition \ref{defwas}.
\end{prop}

\textsc{Proof of Proposition \ref{5lemme_pour_gronwall}.} 
We suppose that $W_2(\rho_1,\rho_2)<+\infty$ otherwise there is nothing to prove.
For $\theta \in [1,2]$, we take $\rho_{\theta}$ that interpolates between $\rho_1$ and $\rho_2$ as in Theorem  \ref{interpol}, and $v_\theta$ the corresponding velocity field.
If we consider, for every $\theta\in[1,2]$, $\Psi_{\theta}$ solution of 
\beq\label{poissontheta}
\Delta \Psi_{\theta}=\rho_{\theta},
\enq
then
$\Psi_{\theta}$ interpolates between $\Psi_1$ and $\Psi_2$.
If we differentiate (\ref{poissontheta}) with respect to $\theta$, we obtain
\be
\Delta \partial_\theta\Psi_\theta &=& \partial_\theta \rho_\theta\\
&=&-\nabla\cdot (\rho_\theta v_\theta).
\en
Note that since $\rho_1, \rho_2$ are bounded, so will be $\rho_\theta$ from Theorem \ref{interpol}, hence using that
$W_2^2(\rho_1,\rho_2) = \int \rho_\theta|v_\theta|^2 < +\infty$, we have $\rho_\theta v_\theta$ bounded in $L^2$.
This implies $\partial_\theta\Psi_\theta \in W^{1,2}$.

We integrate the above equation against $\partial_\theta \Psi_\theta$ to obtain
\be
\int |\nabla\partial_\theta\Psi_\theta|^2 = -\int \rho_\theta v_\theta \cdot \nabla \partial_\theta\Psi_\theta,
\en
and this yields
\be
\|\partial_\theta\nabla\Psi_\theta\|_{L^2(\Rd)} &\leq& \|\rho_\theta v_\theta\|_{L^2(\Rd)}\\
&\leq&\|\rho_\theta\|_{\Linf(\Rd)}^{1/2}W_2(\rho_1,\rho_2).
\en
Using then the second point of Theorem \ref{interpol}, and integrating over $\theta \in [1,2]$, 
this proves Proposition \ref{5lemme_pour_gronwall}.

$\hfill \Box$

\subsection{Lagrangian formulation of the Vlasov-Poisson system}

Given a solution of (VP) with bounded density $\rho$ on $[0,T[$ , we consider for $t\in [0,T[$ the characteristics of equation (\ref{vp1}), that solve the ODE 
\beq
&&\dot X = \xi,\label{ode1}\\
&&\dot \xi = \nabla\Psi(t,X)\label{ode2}.
\enq
Since we assume an $\Linf([0,T[\times\R^3)$ bound on the density $\rho$, the field $\nabla\Psi$ classically satisfies a log-Lispchitz condition:
\beq
&&\forall t\in [0,T[, \ \forall (x,y)\in \R^3\times\R^3, |x-y| \leq \demi, \nonumber \\
&&|\nabla\Psi(t,x) - \nabla\Psi(t,y)| \leq C|x-y|\log\frac{1}{|x-y|},
\enq
where $C$ depends on $\|\rho\|_{\Linf([0,T[\times\R^3)}$.
This condition is enough to define a H\"older continuous flow
\be
\Xi(t,x,v)=(X,\xi)(t,x,v)
\en for the ODE (\ref{ode1},\ref{ode2}), where $(X,\xi)$ is the pair (velocity, position) at time $t$ of the trajectory having (velocity, position) equal to $(x,v)$ at time $0$.

Then we use the following Theorem, proved in \cite{AmbGiSa}:
\begin{theo}
Let $u(t,x)$ be a vector field on $\Rd$. 
Consider the ODE
\be
\dot\gamma(t)=u(t,\gamma(t)),
\en
and the PDE

\be
\dt \mu(t,x)+ \nabla\cdot(\mu(t,x)u(t,x))=0.
\en
Let $B\subset \Rd$ be a Borel set. The following are equivalent:
\begin{enumerate}
\item[(a)] For all $x$ in $B$, there exists a unique solution to the ODE starting at $x$.
\item[(b)] Non negative measure-valued solutions to the PDE with initial data $\mu^0$ concentrated in $B$ are unique.
\end{enumerate}
\end{theo}

From this result, we deduce the following corollary:

\begin{cor}
 The potential $\Psi$ being held fixed, and satisfying
$\Delta \Psi \in \Linf ([0,t]\times\R^3)$, for any $f^0 \in {\cal M}^+(\R^6)$ there exists a unique weak solution to (\ref{vp1}) (i.e. in the sense of (\ref{weakdef}))  with initial datum $f^0$ which is given by
\beq\label{fXi}
f(t) = \Xi(t,\cdot,\cdot)_{\#}f^0,
\enq
where $\Xi=(X,\xi)$ solves (\ref{ode1}, \ref{ode2}).
Note also that we will have 
\beq\label{rhoX}
\rho(t)=X(t,\cdot,\cdot)_{\#} f^0.
\enq
\end{cor}
We remind the reader that the measure $f(t)=\Xi(t, \cdot,\cdot)_{\#}f^0$ is defined by $f(t)(B)=f^0(\Xi^{-1}(t)(B))$ for all Borel subsets $B$ of $\R^6$.

\textsc{Remark.} This corollary does not solve the uiqueness problem, but says only that if we suppress the  coupling between $\Psi$ and $\rho$, there is a unique weak measure-valued solution to the transport equation (\ref{vp1}), that we can represent with the help of characteristics.

\subsection{Final estimate}
Given an initial distribution $f^0(x,\xi)\in {\cal M}^+(\R^6)$ with $\int_{\R^6}  f^0=1$, 
we take two solutions $(f_1,f_2)$ to (VP) with bounded density and initial datum $f^0$. We have $\Delta \Psi= \rho_i, i=1,2$ in the sense of (\ref{2bis}).
 We then consider the associated characteristics $\Xi_1$ and $\Xi_2$, where for $i=1,2$, $\Xi_i=(X_i, \xi_i)(t,x,\xi)$ and $X_i,\xi_i$ solve (\ref{ode1}, \ref{ode2}) with force field $\nabla\Psi_i$.
Note that we will have $f_i(t)=\Xi_i(t)_{\#}f^0, i=1,2$.
We then consider $$Q(t) = \demi \int_{\R^6}  f^0(x,\xi)\left|\Xi_1(t,x,\xi)-\Xi_2(t,x,\xi)\right|^2.$$

\textsc{Remark.}
Notice that $(\Xi_1(t),\Xi_2(t))_{\#}f^0$ is a probability measure on $\R^3\times\R^3$, with marginals $f_1(t)$ and $f_2(t)$, hence by Definition  \ref{defwas}, 
\be
Q(t) \geq \demi W_2^2(f_1(t),f_2(t)).
\en
(We will repeat this argument in Lemma \ref{W2}.) 
This implies in particular that $(Q= 0) \iff (f_1=f_2)$. 
Our proof will rely on an estimate on the Wasserstein distance between $f_1$ and $f_2$, while the proof of \cite{Robert} was obtained by estimating the $H^{-1}$ norm of $f_1-f_2$.

\bigskip
 
Of course $Q(0)=0$, and 
\be
\Dt Q(t) &=& \int_{\R^6}  f^0(x,\xi)(\Xi_1(t,x,\xi)-\Xi_2(t,x,\xi))\cdot \dt (\Xi_1(t,x,\xi)-\Xi_2(t,x,\xi))\\
&= & \int_{\R^6}  f^0(x,\xi)\left[(X_1(t,x,\xi)-X_2(t,x,\xi))\cdot(\xi_1(t,x,\xi)-\xi_2(t,x,\xi)) \right]\\
&+& \int_{\R^6}  f^0(x,\xi)\Big[(\xi_1(t,x,\xi)-\xi_2(t,x,\xi))\cdot \\
&& \hspace{2.5cm} (\nabla\Psi_1(t,X_1(t,x,\xi))-\nabla\Psi_2(t,X_2(t,x,\xi)))\Big].
\en
The second line is bounded by $ Q(t)$, and using Cauchy-Schwartz inequality,  the third line is bounded by
\be
&&Q^{1/2}(t)\left(\int_{\R^6}  f^0(x,\xi)\left|\nabla\Psi_1(t,X_1(t,x,\xi))-\nabla\Psi_2(t,X_2(t,x,\xi))\right|^2\right)^{1/2}\\
&\leq & Q^{1/2}(t)\left(\int_{\R^6}  f^0(x,\xi)\left|\nabla\Psi_2(t,X_1(t,x,\xi))-\nabla\Psi_2(t,X_2(t,x,\xi))\right|^2\right. \\
&& \  \  \  \ \  \  \  \  + \left.\int_{\R^6}  f^0(x,\xi)\left|\nabla\Psi_2(t,X_1(t,x,\xi))-\nabla\Psi_1(t,X_1(t,x,\xi))\right|^2\right)^{1/2} \\
&=& Q^{1/2}(t)\left(T_1(t) +T_2(t)\right)^{1/2}.
\en
Hence we have
\beq\label{intermediaire}
\Dt Q(t) \leq Q(t) +  Q^{1/2}(t)\left(T_1(t) +T_2(t)\right)^{1/2},
\enq
and we will now estimate $T_2$ and then $T_1$.

\bigskip

For $T_2$ we have, using (\ref{rhoX}) and  Proposition \ref{5lemme_pour_gronwall},
\be
T_2(t)&=& \int_{\R^3}  \rho_1(t,x)\left|\nabla\Psi_1(t,x) - \nabla\Psi_2(t,x)\right|^2\\
&\leq & C W_2^2(\rho_1(t),\rho_2(t)),
\en
where $C$ depends on the $\Linf$ norms of $\rho_1,\rho_2$.
Moreover, from the very definition of the Wasserstein distance given in Definition \ref{defwas}, we have the elementary lemma:
\begin{lemme}\label{W2}
Given $X_1, X_2, \rho_1, \rho_2$ as above, we have
\be
W_2(\rho_1(t), \rho_2(t)) \leq \left(\int_{\R^6}f^0(x,\xi)\left|X_1(t,x,\xi)- X_2(t,x,\xi)\right|^2\right)^{\demi}.
\en
\end{lemme}

\textsc{Proof.} The proof follows immediately from Definition \ref{defwas} by noticing that
$\pi = (X_1(t),X_2(t))_{\#} f^0$ is a probability measure on $\R^3\times\R^3$ with marginals $\rho_1(t)$ and $\rho_2(t)$.

$\hfill \Box$

Hence 
\be
W_2^2(\rho_1,\rho_2)& \leq& \int_{\R^6}  f^0(x,\xi)\left|X_2(t,x,\xi)-X_1(t,x,\xi)\right|^2\\ 
&\leq& \int_{\R^6}  f^0(x,\xi)\left|\Xi_1(t,x,\xi)-\Xi_2(t,x,\xi)\right|^2=2Q(t),
\en 
and we conclude that $T_2(t) \leq C Q(t)$ for some $C$ depending on the $\Linf$ bounds on $\rho_1,\rho_2$.

\bigskip

Now, we evaluate $T_1$  by standard arguments,  using the log-Lipschitz regularity of $\nabla\Psi_2$: note first that since $\rho_i, i=1,2$ are bounded in $\Linf$, $\nabla\Psi_i, i=1,2$ are also bounded in $\Linf$, hence for any $C>0$, we can take $T$ small enough such that 
$\|\Xi_i-\Xi_2\|_{\Linf([0,T]\times\R^6)} \leq C$.
Thus we have, for some other $C$ depending on $\|\rho_i\|_{\Linf}, i=1,2$, and as long as $\|\Xi_1-\Xi_2\|_{\Linf}\leq \demi$, 
\be
T_1&=&\int_{\R^6}  f^0(x,\xi)\left|\nabla\Psi_2(t,X_1(t,x,\xi))-\nabla\Psi_2(t,X_2(t,x,\xi))\right|^2\\ 
&\leq& C \int_{\R^6}  f^0(x,\xi) \left(|X_1-X_2|^2 \log^2\frac{1}{|X_1-X_2|}\right)(t,x,\xi) \\
&=&\frac{C}{4}\int_{\R^6}  f^0(x,\xi) \left(|X_1-X_2|^2 \log^2(|X_1-X_2|^2)\right)(t,x,\xi).
\en
Then we use that $x\mapsto x\log^2 x$ is concave for $0\leq x \leq 1/e$,
and we can assume (taking $T$ small enough) that $\|\Xi_1-\Xi_2\|_{\Linf([0,T]\times\R^6)} \leq 1/e$, therefore by Jensen's inequality we have 
\be
&&T_1(t)\\
\leq &&\frac{C}{4}\left[\int_{\R^6}  f^0(x,\xi) |X_1-X_2|^2(t,x,\xi)\right] \, \log^2 \left[\int_{\R^6}  f^0(x,\xi)  |X_1-X_2|^2(t,x,\xi)\right]\\
\leq&& \frac{C}{2} Q(t) \log^2(2Q(t)). 
\en
Combining all these bounds in (\ref{intermediaire}), we obtain that
\be
\Dt Q(t) \leq C Q(t)(1 + \log \frac{1}{Q(t)}),
\en
and we conclude by standard arguments that if $Q(0)=0$, $Q\equiv 0$ on $[0,T[$.

This achieves the proof of Theorem \ref{main1}.

$\hfill \Box$


\bibliography{biblio}
 
{\begin{flushright}{G. Loeper\\EPFL, SB, IMA\\
 10015 Lausanne\\e-mail: {\sf
gregoire.loeper@epfl.ch}}
\end{flushright}}

\end{document}